\documentclass[11pt, a4paper]{amsart}
\usepackage{amsmath, amsthm, amssymb}
\usepackage{txfonts, mathrsfs}
\usepackage[usenames]{color}
\usepackage{psfrag}
\usepackage{graphicx}
\usepackage{enumitem}
\usepackage[british]{babel}
\usepackage{setspace}

\numberwithin{equation}{section}

\usepackage[left=4cm, right=4cm, top=3.4cm, bottom=3cm]{geometry}

\newtheorem{thm}{Theorem}[section]
\newcommand{\bt}{\begin{thm}}
\newcommand{\et}{\end{thm}}

\newtheorem{cor}[thm]{Corollary}
\newcommand{\bc}{\begin{cor}}
\newcommand{\ec}{\end{cor}}

\newtheorem{lem}[thm]{Lemma}
\newcommand{\bl}{\begin{lem}}
\newcommand{\el}{\end{lem}}

\newtheorem{prop}[thm]{Proposition}
\newcommand{\bp}{\begin{prop}}
\newcommand{\ep}{\end{prop}}

\newtheorem{defn}[thm]{Definition}
\newcommand{\bd}{\begin{defn}}
\newcommand{\ed}{\end{defn}}

\newtheorem{rmrk}[thm]{Remark}
\newcommand{\br}{\begin{rmrk}}
\newcommand{\er}{\end{rmrk}}

\newtheorem{quest}[thm]{Question}
\newcommand{\bq}{\begin{quest}}
\newcommand{\eq}{\end{quest}}

\newtheorem{ex}[thm]{Example}
\newcommand{\bex}{\begin{ex}}
\newcommand{\eex}{\end{ex}}

\newcommand{\N}{\mathbb{N}}

\newcommand{\R}{\mathbb{R}}

\newdimen\vintkern\vintkern12pt
\def\vint{-\kern-\vintkern\int}

\newcommand{\hm}{{\mathcal H}}

\newcommand{\diam}{\operatorname{diam}}

\begin{document}
\pagebreak

\bibliographystyle{plain}
\title[]{Canonical parametrizations of metric surfaces of higher topology}


\author{Martin Fitzi}

\address{Kantonsschule Heerbrugg\\ Karl-Völker-Strasse 11\\ 9435 Heerbrugg, Switzerland}
\email{martin.fitzi@ksh.edu}

\author{Damaris Meier}

\address{Department of Mathematics\\ University of Fribourg\\ Chemin du Mus\'ee 23\\ 1700 Fribourg, Switzerland}
\email{damaris.meier@unifr.ch}

\date{\today}

\thanks{Research supported by Swiss National Science Foundation Grant 182423.}

\begin{abstract}
We give an alternate proof to the following generalization of the uniformization theorem by Bonk and Kleiner. Any linearly locally connected and Ahlfors 2-regular closed metric surface is quasisymmetrically equivalent to a model surface of the same topology. Moreover, we show that this is also true for surfaces as above with non-empty boundary and that the corresponding map can be chosen in a canonical way. Our proof is based on a local argument involving the existence of quasisymmetric parametrizations for metric discs as shown in
in a paper of Lytchak and Wenger.
\end{abstract}

\maketitle
\renewcommand{\theequation}{\arabic{section}.\arabic{equation}}
\pagenumbering{arabic}

\section{Introduction and statement of main results}\label{sec:Intro}
\subsection{Introduction}\label{sec:intro-basics}

The classical uniformization theorem states that any oriented Riemannian 2-manifold is conformally diffeomorphic to a model surface of constant curvature. The corresponding map provides a canonical parametrization of said Riemannian surface. An appropriate generalized notion of conformal diffeomorphisms in a non-smooth setting is given by quasisymmetric mappings. A homeomorphism $f\colon X\to Y$ between metric spaces is \emph{quasisymmetric} if there exists a homeomorhism $\eta:[0,\infty)\to[0,\infty)$ such that 
$$d_Y(f(x),f(y))\leq \eta(t)\cdot d_Y(f(x),f(z))$$
for all points $x,y,z\in X$ with $d_X(x,y)\leq t\cdot d_X(x,z)$. The quasisymmetric uniformization problem in the field of analysis on metric spaces then asks under which conditions on a metric space $X$ topologically equivalent to some model space $M$ one may identify $X$ with $M$ via a quasisymmetric homeomorphism. 

A breakthrough result due to Bonk and Kleiner \cite{BK02} asserts that if $X$ is an Ahlfors $2$-regular metric space homeomorphic to the $2$-sphere $S^2$, then there exists a quasi-symmetric homeomorphism between $X$ and $S^2$ if and only if $X$ is linearly locally connected. For definitions of Ahlfors $2$-regularity and linear local connectedness we refer to Section \ref{sec:notations}.

Lytchak and Wenger provide in \cite{LW20} an alternate proof of the theorem of Bonk-Kleiner using a theory of energy and area minimizing discs in metric spaces admitting a quadratic isoperimetric inequality established in \cite{LW15-Plateau} and \cite{LW-intrinsic}. The aim of this paper is to use the existence result in \cite{LW20} locally to obtain canonical parametrizations of metric surfaces of higher topology with possibly non-empty boundary.

Let $X$ be a metric space homeomorphic to a smooth surface $M$. Here, a \emph{smooth surface} referes to a smooth compact oriented and connected Riemannian $2$-manifold with possibly non-empty boundary. Define $\Lambda(M,X)$ to be the family of Newton-Sobolev maps $u\in N^{1,2}(M,X)$ such that $u$ is a uniform limit of homeomorphisms from $M$ to $X$ and let $E^2_+(u,g)$ be the \emph{Reshetnyak energy} of a map $u\in N^{1,2}(M,X)$ with respect to the Riemannian metric $g$; for definitions see Section~\ref{sec:introSobolev}. Our main result is the following version of \cite[Theorem~1.1]{LW20} for metric surfaces of higher topology. Note that the definition of $\Lambda(M,X)$ is different from \cite{LW20}.

\bt 
\label{thm:main}
Let $X$ be a metric space which is Ahlfors $2$-regular, linearly locally connected and homeomorphic to a smooth surface $M$. 
Then, there exist a map $u\in \Lambda(M,X)$ and a Riemannian metric $g$ on $M$ such that
$$E_+^2(u,g)=\inf\{E_+^2(v,h):v\in \Lambda(M,X),\, h\text{ a smooth Riemannian metric on M}\}.$$
Any such $u$ is a quasisymmetric homeomorphism from $M$ to $X$ and the pair $(u,g)$ is uniquely determined up to a conformal diffeomorphism $\varphi\colon (M, g)\to (M,h)$.
\et 

Moreover, the metric $g$ can be chosen to be of constant sectional curvature $-1$, 0 or 1 and such that $\partial M$ is geodesic (if non-empty).

The theorem of Bonk-Kleiner has been extended for example in \cite{Wildrick08}, \cite{MW13} and \cite{Wildrick10}. In particular, there are analogous quasisymmetric uniformization results for Ahlfors $2$-regular (orientable) closed surfaces, see \cite{GW18} and \cite{Iko19}. Theorem~\ref{thm:main} is a strengthening of these results, since it states the existence of \emph{canonical} quasisymmetric homeomorphisms, even for such surfaces with boundary. Hence we obtain the following generalization of the result of Bonk-Kleiner.
 
\bc
\label{cor:genOfBonkKleiner}
Let $X$ be an Ahlfors $2$-regular metric space homeomorphic to a smooth surface $M$ with possibly non-empty boundary. Then, $X$ is quasisymmetrically equivalent to $M$ if and only if $X$ is linearly locally connected.
\ec

\subsection{Elements of proof} We briefly sketch some of the arguments needed for proving Theorem~\ref{thm:main}. For arbitrary $M$ and $X$ as in the paragraph before the theorem, the set $\Lambda(M,X)$ can be empty. A crucial step in this work is to show the existence of a map $u\in\Lambda(M,X)$ in the setting of Theorem~\ref{thm:main}.

\bp
\label{prop:lambdaNotEmpty}
Let $M$ be a smooth surface and $(X,d)$ a metric space which is Ahlfors $2$-regular, linearly locally connected and homeomorphic to $M$. Then the family $\Lambda(M,X)$ is non-empty and contains a quasisymmetric homeomorphism.
\ep

The proposition follows by a dissection of $M$ and $X$ into appropriate disc-type subdomains, consequently applying \cite[Theorem~6.1]{LW20} yielding quasisymmetric parametrizations for each subdomain in $X$ and finally gluing all these mappings together in order to obtain a global quasisymmetric homeomorphism $M\to X$. Note that Proposition~\ref{prop:lambdaNotEmpty} already establishes Corollary~\ref{cor:genOfBonkKleiner}.

The map $u$ provided by Proposition~\ref{prop:lambdaNotEmpty} is not necessarily canonical, i.e. possibly not of minimal energy. In order to find an energy minimizer in $\Lambda(M,X)$, we will use similar arguments as in the proofs of \cite[Theorem~6.1]{LW20} and \cite[Theorem~8.2]{FW-Plateau-Douglas}. In particular, we need to ensure that a family of mappings in $\Lambda(M,X)$ of uniformly bounded energies is equicontinuous. 

The paper is structured as follows. In Section~2 we provide necessary definitions and some results on Newtonian Sobolev spaces that will be of use later on. Section~3 is devoted to the proof of Proposition~\ref{prop:lambdaNotEmpty}. In Section~4 we will show equicontinuity of energy bounded almost homeomorphisms. And finally, the proof of Theorem~\ref{thm:main} is given in Section~5.

\bigskip
{\bf Acknowledgments:} We wish to thank our PhD advisor Stefan Wenger for his great support and numerous discussions on this topic.

\section{Preliminaries} 

\subsection{Basic definitions and notations}
\label{sec:notations}
Let $(X,d)$ be a metric space. The \emph{open ball} in $X$ of radius $r>0$ centered at a point $x\in X$ is denoted by $B_X(x,r)$ or simply $B(x,r)$. Consider the Euclidean space $(\R^2,|\cdot|)$, where $|\cdot|$ is the Euclidean norm. The \emph{open and closed unit discs} in $\R^2$ are given by
$$D:=\{z\in\R^2:|z|<1\},\qquad \overline{D}:=\{z\in\R^2:|z|\leq1\}.$$
An open set $\Omega\subset X$ homeomorphic to the unit disc $D$ is a \emph{Jordan domain in $X$} if its completion $\overline{\Omega}\subset X$ is homeomorphic to $\overline{D}$. A \emph{Jordan curve} in $X$ is a subset of $X$ homeomorphic to $S^1$ and it is called \emph{chord-arc} if it is biLipschitz equivalent to $S^1$. The \emph{image} of a curve $c$ in $X$ is denoted by $|c|$ and the \emph{length} by $\ell_X(c)$ or $\ell(c)$. A curve $c\colon[a,b]\to X$ is called \emph{geodesic} if $\ell(c)=d(c(a),c(b))$. A metric space $(X,d)$ is \emph{geodesic} if every pair of points in $X$ can be joined by a geodesic. 

A \emph{metric surface} $X$ is a metric space homeomorphic to a smooth surface $M$. We write $\partial X$ and $\partial M$ for the boundaries of $X$ and $M$, respectively.
 
For $s\geq 0$, we denote the \emph{$s$-dimensional Hausdorff
measure} of a set $A\subset X$ by $\hm_X^s(A)$ or simply $\hm^s(A)$. The normalizing constant is chosen in such a way that if $X$ is the Euclidean space $\R^n$, the Lebesgue measure agrees with $\hm^n_{X}$. If $(M,g)$ is a Riemannian manifold of dimension $n$ then the $n$-dimensional Hausdorff measure $\hm_g^n:=\hm_{(M,g)}^n$ on $(M,g)$ coincides with the Riemannian volume.

Let $g$ be a smooth Riemannian metric on a smooth surface $M$ such that the boundary of $M$ is geodesic with respect to $g$. We call the metric $g$ \emph{hyperbolic} if it is of constant sectional curvature $-1$, and \emph{flat} if it has vanishing sectional curvature as well as an associated Riemannian 2-volume satisfying $\hm^2_g(M)=1$.

\bd
A metric space $X$ is said to be \emph{Ahlfors $2$-regular} if there exists $K>0$ such that for all $x\in X$ and $0<r<\diam X$, we have
$$K^{-1}\cdot r^2\leq \hm^2_X(B(x,r))\leq K\cdot r^2.$$
\ed

\bd We say that a metric space $X$ is \emph{linearly locally connected (LLC)} if there exists a constant $\lambda\geq1$ such that for all $x\in X$ and $r>0$, every pair of distinct points in $B(x,r)$ can be connected by a continuum in $B(x,\lambda r)$ and every pair of distinct points in $X\setminus B(x,r)$ can be connected by a continuum in $X\setminus B(x,r/\lambda)$.  
\ed

\subsection{Metric space valued Sobolev maps}
\label{sec:introSobolev}
In this subsection we give a brief overview over some basic concepts used in the theory of metric space valued Sobolev maps based on upper gradients. Note that several other equivalent definitions of Sobolev spaces exist. For more details consider e.g. \cite{HKST15}. 

Let $(X,d)$ be a complete metric space and $M$ a smooth compact $2$-dimensional manifold, possibly with non-empty boundary. Fix a Riemannian metric $g$ on $M$ and consider a domain $\Omega\subset M$. Let $u\colon\Omega\to X$ be a map and $\rho\colon\Omega\to [0,\infty]$ a Borel function. Then, $\rho$ is called \emph{(weak) upper gradient of $u$ with respect to $g$} if
\begin{align}
    \label{ineq:UpperGradient}
    d(u(\gamma(a)),u(\gamma(b)))\leq\int_{\gamma}\rho(s)\;ds
\end{align}
for (almost) every rectifiable curve $\gamma\colon[a,b]\to\Omega$. A weak upper gradient $\rho$ of $u$ is said to be \emph{minimal} if $\rho\in L^2(\Omega)$ and for every weak upper gradient $\rho'$ of $u$ in $L^2(\Omega)$ it holds that $\rho\leq\rho'$ almost everywhere on $\Omega$. 
Denote by $L^2(\Omega,X)$ the family of measurable essentially separably valued maps $u\colon\Omega\to X$ such that the distance function $u_x(z):=d(u(z),x)$ is in the space $L^2(\Omega)$ of $2$-integrable functions for some and hence any $x\in X$. A sequence $(u_k)\subset L^2(\Omega,X)$ is said to \emph{converge in $L^2(\Omega,X)$} to a map $u\in L^2(\Omega,X)$ if
$$\int_{\Omega}d^2(u_k(z),u(z))\;d\hm^2_g(z)\to0$$
as $k$ tends to infinity. The \emph{(Newton-)Sobolev space} $N^{1,2}(\Omega,X)$ is the collection of maps $u\in L^2(\Omega,X)$ such that $u$ has a weak upper gradient in $L^2(\Omega)$. Every such $u$ has a minimal weak upper gradient denoted by $\rho_u$, which is unique up to sets of measure zero (see e.g. \cite[Theorem~6.3.20]{HKST15}). Note also that the definition of $N^{1,2}(\Omega,X)$ is independant of the chosen metric $g$ on $M$.

\bd
The \emph{Reshetnyak energy} of a map $u\in N^{1,2}(\Omega,X)$ with respect to $g$ is defined by
$$E^2_+(u,g):=\int_{\Omega}|\rho_u(z)|^2\;d\hm_g^2(z).$$
\ed

This definition of energy agrees with the one given in \cite[Definition 2.2]{FW-Plateau-Douglas}; in particular, $E_+^2$ is invariant under precompositions with conformal diffeomorphisms.

\section{Noncanonical quasisymmetric parametrizations}\label{sec:gluing}

\subsection{Decompositions of metric Y-pieces and cylinders into Jordan domains}

We begin this section by introducing some terminology. A \textit{cylinder} and \textit{Y-piece} are connected topological surfaces of genus 0  with two and three boundary components, respectively. Furthermore, we refer to a metric space homeomorphic to a cylinder or a Y-piece as a \textit{metric cylinder} or a \textit{metric Y-piece}, respectively.

\bl
\label{lem:JordanDecomp}
Let $X$ be a geodesic metric surface and $\Sigma\subset X$ a metric cylinder or metric Y-piece such that each connected component of $\partial\Sigma$ can be parametrized by a piecewise geodesic chord-arc curve. Then there exist Jordan domains $J_1,J_2\subset \Sigma$ with
\begin{enumerate}
    \item $\Sigma=\overline{J_1}\cup\overline{J_2}$,
    \item $J_1\cap J_2=\emptyset$,
    \item $J_1,J_2$ are both bounded by a biLipschitz curve.
\end{enumerate}
\el

\begin{proof}
We give a proof for $\Sigma$ being a metric Y-piece, the case of a metric cylinder only needing minor adaptations in the following arguments. Denote by $\eta_i\colon S^1\to \partial\Sigma$ the piecewise geodesic biLipschitz curves parametrizing the three components of $\partial \Sigma$. Choose three disjoint injective curves $\gamma_i$ in $\Sigma$, each one connecting two boundary components such that $\Sigma$ is separated into two Jordan domains when cutting along these curves. By \cite[Lemma 4.2]{LW20} and its proof, we may assume that each $\gamma_i$ is biLipschitz and piecewise geodesic. Denote the endpoints of $\gamma_i$ by $a_i^1$, $a_i^2$.

Choose $\varepsilon>0$ so small that the balls $B(a_i^j, 2\varepsilon)$ are disjoint. We modify $\gamma_i$ within $B(a_i^j, 2\varepsilon)$ with the following procedure. Without loss of generality assume $a_1^1\in |\eta_1|$. Choose a point $x_1\in B(a_1^1,\varepsilon)\cap |\gamma_1|$ distinct from $a_1^1$ and let $y_1\in|\eta_1|$ be such that
\begin{align}
    \label{eq:dist}
    d(x_1,y_1)=d(x_1,|\eta_1|),
\end{align}
where $d$ denotes the metric on $X$. Let $c_1\colon I\to  \Sigma$ be a geodesic segment connecting $x_1$ with $y_1$. Thus, $|c_1|\subset B(a_1^1, 2\varepsilon)$. Then consider the concatenation of $c_1$ with one of the subcurves of $\eta_1$ emanating from $y_1$. Let $s>0$ be such that the following holds. Subcurves of $\eta_1$ and $c_1$ with common endpoint $y_1$ can be reparametrized by arc-length on $[-s,0]$ and $[0,s]$, respectively, such that $\eta_1(0)=c_1(0)=y_1$. Denote this concatenation defined on $[-s,s]$ by $\eta$. Equality (\ref{eq:dist}) implies that for $r\in[0,s]$
$$d(\eta_1(-r), c_1(r))\geq r.$$
It follows from the proof of \cite[Lemma 4.2]{LW20} that $\eta$ is a biLipschitz curve. Redefine $\gamma_1$ by replacing the subcurve from $x_1$ to $a_1^1$ by $c_1$. Analogously, construct segments $c_2,\dots,c_6$ in the vicinities of the other $a_i^j$ and modify every $\gamma_i$ near its endpoints in this way. By choosing appropriate subcurves, we have that all redefined $\gamma_i$ are still injective. Moreover, the proof of \cite[Lemma 4.2]{LW20} shows that if $\gamma_i$ is not biLipschitz at a vertex in the interior of the curve, we can change it in an arbitrarily small ball around this vertex to obtain a global biLipschitz curve.

Finally, $\Sigma$ is separated into Jordan domains $J_1$ and $J_2$ by cutting along redefined $\gamma_i$. Moreover, the boundaries $\partial J_1$ and $\partial J_2$ are parametrized by biLipschitz concatenations of the redefined $\gamma_i$ with respective subcurves of $\eta_j$.
\end{proof}

The following lemma will be useful in the proof of Proposition~\ref{prop:QSExtensionCylinder}. A \textit{metric disc} is a metric space homeomorphic to the closed unit disc $\overline{D}$.

\bl\label{lem:subsetLLCRegular}
Let $X$ be an Ahlfors $2$-regular and LLC metric surface. Consider a chord-arc curve $\gamma$ in $X$. Let $\Sigma\subset X$ be a metric disc bounded by $\gamma$ or a metric cylinder with boundary components $\partial\Sigma^1=|\gamma|$ and $\partial\Sigma^2\subset\partial X$. Then $\Sigma$ equipped with the subspace metric is Ahlfors $2$-regular and LLC.
\el

The lemma can be shown readily by using the LLC-property of $X$ and replacing parts of the continua which lie in $X\setminus\Sigma$ with appropriate subcurves of the biLipschitz curve at the boundary $\gamma$ in order to obtain desired continua in $\Sigma$. Compare also to the proof of \cite[Proposition 6.4]{LW20}.

\subsection{Parametrizations of boundary cylinders}
The aim of this section is to establish the following extension result for cylindrical surfaces which is needed later in the proof of Proposition~\ref{prop:lambdaNotEmpty}.

\bp\label{prop:QSExtensionCylinder}
Let $Z$ be a smooth cylinder and $\partial Z^1\subset\partial Z$ a boundary component. Let $\Sigma$ be an Ahlfors $2$-regular and LLC metric cylinder and $\partial \Sigma^1\subset\partial \Sigma$ a boundary component. Assume furthermore that there exists a biLipschitz homeomorphism $f\colon \partial Z^1\to\partial \Sigma^1$. Then $f$ extends to a quasisymmetric homeomorphism $\overline{f}\in\Lambda(Z,\Sigma)$.
\ep

As a first step in the proof of Proposition~\ref{prop:QSExtensionCylinder}, we will perform a gluing of the metric cylinder $\Sigma$ with the closed unit disc $\overline{D}$ along corresponding boundary components. We now introduce some notation and needed results concerning this gluing method.

Let $(X,d_X)$ and $(Y,d_Y)$ be two compact metric surfaces with non-empty boundary and let $\partial X^j\subset\partial X$, $\partial Y^k\subset\partial Y$ be two boundary components. Assume $\gamma\colon \partial X^j\to\partial Y^k$ is a biLipschitz homeomorphism and define the quotient
$$\widehat{XY}:=(X\sqcup Y)/\sim,$$
where $x\sim y$ for $x\in X$, $y\in Y$ if $y=\gamma(x)$. Equip $\widehat{XY}$ with the quotient metric, denoted $\hat{d}$ (see e.g. \cite[Definition~3.1.12]{BBI01}). Consider $X$ and $Y$ as subsets of $\widehat{XY}$ and set $X\cap Y:=\{[x]:x\in\partial X^j\}$. It follows immediately that the identity maps $(X,d_X)\to(X,\hat{d}|_{X\times X})$ and $(Y,d_Y)\to(Y,\hat{d}|_{Y\times Y})$ are 1-Lipschitz. The next lemma is a consequence of the compactness of $X\cap Y$ and the biLipschitz property of $\gamma$.

\bl\label{lem:identity-L-Lip}
The identity maps $\mathrm{id}_X\colon(X,\hat{d}|_{X\times X})\to(X,d_X)$ and $\mathrm{id}_Y\colon(Y,\hat{d}|_{Y\times Y})\to(Y,d_Y)$ are $L$-Lipschitz, where $L\geq 1$ denotes the biLipschitz constant of $\gamma$. In particular, the restrictions $\hat{d}|_{X\times X}$ and $\hat{d}|_{Y\times Y}$ are $L$-biLipschitz equivalent to $d_X$ and $d_Y$. 
\el

Moreover, we have the following geometric property of the space $(\widehat{XY}, \hat{d})$.

\bl 
\label{lem:QuotientRegularLLC}
If $(X,d_X)$ and $(Y,d_Y)$ are Ahlfors $2$-regular and LLC, then so is $(\widehat{XY}, \hat{d})$.
\el

The proof of Lemma~\ref{lem:QuotientRegularLLC} can be found in the appendix.
Another important result used in the proof of Proposition~\ref{prop:QSExtensionCylinder} (and the proof of Proposition~\ref{prop:lambdaNotEmpty}) is the following variant of \cite[Theorem 6.1]{LW20}. 

\bt
\label{thm:QSExtensionBiLipBoundary}
Let $X$ be an Ahlfors $2$-regular geodesic metric space homeomorphic to a 2-dimensional manifold. Let $J\subset X$ be a Jordan domain with $\ell(\partial J)<\infty$ and such that $\overline{J}$ is LLC. Then any quasisymmetric homeomorphism $f\colon S^1\to\partial J$ extends to a quasisymmetric homeomorphism $\overline{f}\in\Lambda(\overline{D},\overline{J})$. 
\et

Note that the parametrization of $u$ at $\partial \Omega$ can be prescribed by precomposing with a suitable quasisymmetric homeomorphism of $\overline{\Omega}$ obtained via the extension result in \cite{BA56}. See also the proof of \cite[Proposition 6.4]{LW20} or compare to the similar use of a variant of the aforementioned extension result for mappings between boundary circles of annuli-type surfaces, \cite[Theorem 3.14]{TV81}, at the end of the next proof.

\begin{proof}[Proof of Proposition~\ref{prop:QSExtensionCylinder}]
Consider the quotient space
$$\widehat{\Sigma D}:=(\Sigma\sqcup \overline{D})/\sim$$
defined as above for some biLipschitz homeomorphism $\partial \Sigma^1\to \partial D=S^1$ and equipped again with the quotient metric $\widehat{d}$. By Lemma~\ref{lem:QuotientRegularLLC}, the metric disc $(\widehat{\Sigma D},\widehat{d})$ is Ahlfors $2$-regular and LLC and hence quasiconvex (see \cite[Theorem B.6]{Sem96}). Therefore, the space $\widehat{\Sigma D}$ is geodesic up to a biLipschitz change of metric and by Theorem~\ref{thm:QSExtensionBiLipBoundary} there exists a quasisymmetric homeomorphism $v\in\Lambda(\overline{D}, \widehat{\Sigma D})$.
Consider $\widehat{D}=D$ as a subset of $\widehat{\Sigma D}$ and define
$$\Omega:=\overline{D}\backslash v^{-1}(\widehat{D}). $$
By the annulus conjecture (see \cite[Theorem 3.12]{TV81}) there exists a quasisymmetric homeomorphism $g\colon A \to \Omega$, where
$$A:=\{p\in\R^2: 1/2\leq |p|\leq 1\}\subset \overline{D}$$
denotes the standard annulus equipped with the Euclidean metric. Without loss of generality, we may assume that $g$ maps the unit circle onto $\partial (v^{-1}(\widehat{D}))$. Let $\varphi\colon Z\to A$ be a biLipschitz homeomorphism with $\varphi(\partial Z^1)=S^1$. Then, the mapping $u\in N^{1,2}(Z,\Sigma)$ defined by $u:= \mathrm{id}_\Sigma \circ v \circ g\circ\varphi$ is a quasisymmetric homeomorphism with $u(\partial Z^1)=\partial\Sigma^1$. Moreover, the composition
$$h:=\varphi\circ u^{-1}\circ f\circ\varphi^{-1}|_{S^1}\colon S^1\to S^1$$
is a quasisymmetric homeomorphism, which we may assume to be orientation-preser\-ving. By \cite[Theorem~3.14]{TV81}, the map $h$ extends to a quasisymmetric homeomorphism $\overline{h}\colon\overline{D}\to\overline{D}$ such that $\overline{h}$ restricted to the ball $B(0,1/2)$ is the identity map. Hence 
$$\overline{f}:=u\circ \varphi^{-1}\circ\overline{h}\circ \varphi$$
is a desired quasisymmetric homeomorphism from $Z$ to $\Sigma$ with $\overline{f}|_{\partial Z^1}=f$.
\end{proof}

\subsection{Noncanonical quasisymmetric parametrizations} 

Using the extension result established in the previous subsection, we may obtain Proposition~\ref{prop:lambdaNotEmpty} mentioned in the introduction. 

\begin{proof}[Proof of Proposition \ref{prop:lambdaNotEmpty}]
The cases where $M$ is a disc or a sphere follow from \cite[Theorem 6.1]{LW20} and \cite[Proposition 6.4]{LW20}.

Depending on its topology, endow $M$ with a hyperbolic or flat Riemannian metric. Furthermore, we may assume that $X$ is geodesic up to a biLipschitz change of metric (see \cite[Theorem B.6]{Sem96}). Let $h\colon M\to X$ be a homeomorphism.

We first give a proof in the special case when $X$ has either empty boundary or else is bounded by piecewise geodesic chord-arc curves. Choose a collection of simple closed geodesics $\{\gamma_i\colon S^1\to M\}$ decomposing $M$ into smooth Y-pieces or cylinders $M_k$, respectively. Using \cite[Lemma 4.2]{LW20}, we may partition $X$ into Y-pieces/cylinders $X_k$ such that each $X_k$ is homotopic to $h(M_k)\subset X$ and bounded by piecewise geodesic chord-arc curves. We then further decompose $M_k$ and $X_k$ into Jordan domains: if $M_k$ is a Y-piece, then it is a standard result from hyperbolic geometry that $M_k$ is isometric to the partial gluing of the boundary of two copies $\Omega_{k,1}, \Omega_{k,2}$ of a right-angled hexagon in $\mathbb{H}$, see e.g. \cite[Proposition 3.1.5]{Bus10}. If $M_k$ is of cylindrical type, then a similar decomposition into isometric rectangles in the Euclidean plane, again denoted $\Omega_{k,1}$ and $\Omega_{k,2}$, is possible. Note that in either case $\Omega_{k,1}$ and $\Omega_{k,2}$ are biLipschitz equivalent to the closed unit disc $\overline{D}$. In $X$ we decompose each $X_k$ into Jordan domains $J_{k,1}$, $J_{k,2}$ as in Lemma~\ref{lem:JordanDecomp}. After possibly inverting the notation of $J_{k,1}$ and $J_{k,2}$, let $$f\colon \bigcup_{\substack{j=1,2\\ k}} \partial \Omega_{k,j}\to \bigcup_{\substack{j=1,2\\ k}}\partial J_{k,j}$$ be a biLipschitz homeomorphism satisfying $f(\partial \Omega_{k,j})=\partial J_{k,j}$ for each $j$, $k$. By Lemma~\ref{lem:subsetLLCRegular} and Theorem~\ref{thm:QSExtensionBiLipBoundary}, there exists for each $k$ a quasisymmetric homeomorphism $g_{k,j}\colon \overline{\Omega_{k,j}}\to \overline{J_{k,j}}$ with $g_{k,j}|_{\partial \Omega_{k,j}}=f|_{\partial \Omega_{k,j}}$.
The map $u\colon M\to X$ agreeing with $g_{k,j}$ on $\Omega_{k,j}$ satisfies the hypotheses of the quasisymmetric gluing theorem \cite[Theorem~3.1]{AKT05} as each $\Omega_{k,j}$ is bounded and has biLipschitz boundary and every $g_{k,j}$ is a quasisymmetric homeomorphism. Therefore, the map $u$ itself is a quasisymmetric homeomorphism. This shows the proposition in the special case.

We now turn to the general case, where $X$ might be bounded by curves of unknown regularity. For each boundary component $\partial X^i$, define a piecewise geodesic biLipschitz curve $c_i\colon S^1\to X$ which is homotopic to an oriented parametrization of $\partial X^i$, but disjoint from it. Furthermore, we may assume that the curves $\{c_i\}$ are all pairwise disjoint. Let $\Sigma_i\subset X$ be the metric cylinder bounded by $c_i(S^1)$ and $\partial X^i$, and let $\Sigma\subset X$ be the subsurface bounded by $\bigcup_i c_i(S^1)$. Note that $\Sigma$ is homeomorphic to $X$. The first part of the proof then shows that there exists a quasisymmetric homeomorphism $u\colon M\to \Sigma$. Then embed $M$ smoothly into a surface $\Tilde{M}$ such that for each $i$, there exists exactly one boundary component $\partial\Tilde{M}^i$ which together with $\partial Z_i^1:=u^{-1}(c_i(S^1))\subset\partial M$ bounds a smooth cylinder $Z_i\subset\Tilde{M}$. Finally, use Lemma~\ref{lem:subsetLLCRegular} and Proposition~\ref{prop:QSExtensionCylinder} to obtain quasisymmetric extensions $u_i\colon Z_i\to \Sigma_i$ of $u|_{\partial Z_i^1}$. Once again, the gluing result \cite[Theorem~3.1]{AKT05} ensures that the map $u\colon \Tilde{M}\to X$ agreeing with $u$ on $M$ and with $u_i$ on $Z_i$ is a quasisymmetric homeomorphism. The proof of the proposition is complete.
\end{proof}

\section{Equicontinuity of energy bounded almost homeomorphisms}

The map provided by Proposition~\ref{prop:lambdaNotEmpty} does not need to be canonical, i.e. of minimal energy. In order to obtain such a parametrization in Section~\ref{sec:proofOfMainThm}, we will apply a direct variational method for which we need to know equicontinuity of a given energy-minimizing sequence of parametrizations. More explicitly, we prove the following statement in this section.

\bp
\label{prop:equicontinuity}
Let $M$ be a smooth surface endowed with a Riemannian metric $g$ and which is neither of disc- nor of sphere-type. Let $X$ be a metric surface homeomorphic to $M$ and such that $\partial X$ is rectifiable. Then the family
$$\mathcal{F}:=\{v\in\Lambda(M,X):E_+^2(v,g)\leq K\}$$
is equicontinuous.
\ep

In order to show Proposition \ref{prop:equicontinuity}, we need the following elementary lemma. Its proof is left to the reader.

\bl 
\label{lem:diamJordanDomain}
Let $X$ be a metric surface which is not of sphere-type. Then for every $\varepsilon>0$ there exists $\rho>0$ such that the following holds. Every embedding $u\colon\overline{D}\to X$ with $\diam (u(S^1))<\rho$ satisfies $\diam(u(\overline{D}))<\varepsilon$.
\el

By continuity, the statement holds for any uniform limit of embeddings from $\overline{D}$ to $X$.

\begin{proof}[Proof of Proposition \ref{prop:equicontinuity}]
Let $\varepsilon>0$ and define $$\eta:=\inf\{\ell(c)\mid c\colon S^1\to X\text{ is a non-contractible curve in }X\}>0.$$ By Lemma \ref{lem:diamJordanDomain}, there exists $0<\rho<\min\{\varepsilon, \eta\}$ such that for any uniform limit of embeddings $u\colon\overline{D}\to X$ with $\diam(u(S^1))<\rho$ there holds $\diam(u(\overline{D}))<\varepsilon$. Similarly, there exists $0<\rho'<\rho/2$ such that the following is true. If $x,x'\in\partial X$ satisfy $d(x,x')<\rho'$, then they lie on the same component $\partial X^i\subset\partial X$ and the shorter of the two subcurves of $\partial X^i$ connecting $x$ and $x'$ has length at most $\rho/2$. Since $M$ is compact, there exists $0<\delta<1$ so small that
$$\pi\cdot\left(\frac{8K}{|\log(\delta)|}\right)^{1/2}<\rho'$$
and such that every point $p\in M$ is contained in a neighbourhood in $M$ which is the image of the set $B:=B_{\R^2}(q,\sqrt{\delta})\cap\overline{D}$ under a map $\psi$ that is 2-biLipschitz and takes the point $q\in[0,1]\subset\overline{D}$ to $p$, where $q$ is chosen to be $1$ if $p\in\partial M$ and $0$ whenever the distance between $p$ and $\partial M$ is big enough. In particular, if the set $B_{\R^2}(q,\sqrt{\delta})\cap S^1$ is not empty, then it is mapped onto a subcurve of $\partial M$.

Fix $p\in M$ and $v\in\mathcal{F}$. By the Courant-Lebesgue Lemma (see e.g. \cite[Lemma 7.3]{LW15-Plateau}) there exists $r\in (\delta, \sqrt{\delta})$ such that
$$\ell(v\circ\psi\circ\gamma_r)\leq\pi\cdot\left(\frac{2E_+^2(v\circ\psi)}{|\log(\delta)|}\right)^{1/2}\leq\pi\cdot\left(\frac{8E_+^2(v)}{|\log(\delta)|}\right)^{1/2} <\rho',$$
where $\gamma_r$ is an arc-length parametrization of $\{z\in B:|z-q|=r\}$. 

Consider the set $A:=\{z\in B:|z-q|<r\}$. It holds that $B_M(p,\delta/2 )\subset \psi(A)$ and $\overline{A}$ is biLipschitz equivalent to $\overline{D}$ with constant only depending on $r$. If $\psi(A)$ does not intersect $\partial M$, by applying Lemma \ref{lem:diamJordanDomain}, we can conclude $\diam(v(\psi(A)))<\varepsilon$ and therefore $v(B_M(p,\delta/2 ))\subset B_X(v(p),\varepsilon)$.

If $\psi(A)\cap\partial M$ is not empty, then $\psi(A)$ is bounded by $\psi\circ\gamma_r$ and a subarc of $\partial M^i$, denoted $\alpha_r$. The endpoints $a_r,b_r\in \partial M^i$ of $\psi\circ\gamma_r$ satisfy $d(v(a_r),v(b_r))<\rho'<\rho/2$. Thus, $v(a_r)$ and $v(b_r)$ lie on the same boundary component $\partial X^i\subset\partial X$ and the shorter subcurve of $\partial X^i$ connecting $v(a_r)$ and $v(b_r)$ has length at most $\rho/2<\eta/2$. This segment corresponds to the curve $v\circ\alpha_r$. Indeed otherwise, the concatenation of $v\circ\psi\circ\gamma_r$ with $v|_{\partial M^i\setminus \alpha_r}$ would yield a non-contractible closed curve in $X$ of length strictly less than $\eta$, which is impossible. Again by applying Lemma~\ref{lem:diamJordanDomain} we obtain $v(B_M(p,\delta/2 ))\subset v(\psi(A))\subset B_X(v(p),\varepsilon)$. Since the choice of $\delta$ was independant of $v$ and of $p$, this proves equicontinuity of $\mathcal{F}$.
\end{proof}

\section{Proof of Main Theorem}\label{sec:proofOfMainThm}
We finally turn to the proof of Theorem~\ref{thm:main}. First however, we introduce some notation.
Define the family
$$\Lambda_{\text{metr}}(M,X):=\{(v,h):v\in\Lambda(M,X),\, h\text{ a smooth Riemannian metric on M}\}.$$
An \textit{energy minimizing sequence} in $\Lambda_{\text{metr}}(M,X)$ is a sequence of pairs $(u_n,g_n)\in\Lambda_{\text{metr}}(M,X)$ satisfying 
$$E_+^2(u_n,g_n)\to\inf\{E_+^2(v,h):(v,h)\in\Lambda_{\text{metr}}(M,X)\}$$
as $n$ tends to infinity.

\begin{proof}[Proof of Theorem \ref{thm:main}]
The proofs in the cases where $M$ is of disc- or sphere-type follow from \cite{LW20}, and we therefore assume that that $M$ is of higher topological type. In a first step, we show the existence of an energy minimizing pair in $\Lambda_{\text{metr}}$. By Proposition \ref{prop:lambdaNotEmpty}, the set $\Lambda(M,X)$ is not empty. Therefore, we are able to consider an energy minimizing sequence $(u_n,g_n)$ in $\Lambda_{\text{metr}}(M,X)$. We lose no generality in assuming that the metrics $g_n$ are all hyperbolic respectively flat. Observe that each $u_n$, being a uniform limit of homeomorphisms, satisfies the condition of cohesion for some $\eta>0$ in the sense of \cite[Definition 8.1]{FW-Plateau-Douglas}. Thus by \cite[Proposition 8.4]{FW-Plateau-Douglas} there exists $\varepsilon>0$ depending only on $\eta$ and $K:=\sup_{n\in\N}E_+^2(u_n,g_n)$ such that for every $n$ the relative systole of $(M,g_n)$ (see \cite[Definition 3.1]{FW-Plateau-Douglas}) is bounded from below by $\varepsilon$. Then, there exist diffeomorphisms $\varphi_n\colon M\to M$ such that a subsequence of $(\varphi_n^*g_n)$ converges smoothly to a hyperbolic respectively flat metric $g$ on $M$ (see \cite[Theorem 4.4.1]{DHT10} if $M$ is a closed surface; and e.g. \cite[Theorem 3.3]{FW-Plateau-Douglas} if $M$ has non-empty boundary and admits hyperbolic metrics). Set $v_n:=u_n\circ \varphi_n$. The convergence above implies that
$$\lim_{n\to\infty}E_+^2(v_n,g)=\lim_{n\to\infty}E_+^2(u_n,g_n).$$
Thus, the sequence $(v_n,g)$ is energy minimizing in $\Lambda_{\text{metr}}(M,X)$. Now by Proposition~\ref{prop:equicontinuity}, the sequence $(v_n)$ is equicontinuous and the Arzelà-Ascoli theorem implies that a subsequence of $(v_n)$ converges uniformly to some continuous map $u\colon M\to X$. It follows that $u$ is in $N^{1,2}(M,X)$ (compare to \cite[Theorem~1.6.1]{KS93}) as well as a uniform limit of homeomorphisms, hence $u\in\Lambda(M,X)$. By the lower semicontinuity of $E_+^2(\cdot)$ it follows that the pair $(u,g)$ is an energy minimizer in $\Lambda_{\text{metr}}(M,X)$.

We now show that any energy minimizing pair $(u,g)$ in $\Lambda_{\text{metr}}(M,X)$ is a quasisymmetric homeomorphism. As a uniform limit of homeomorphisms, the map $u$ is continuous, monotone and surjective. Furthermore, by \cite[Theorem 4.2]{FW-Plateau-Douglas}, the map $u$ is infinitesimally isotropic and hence infinitesimally $\sqrt{2}$-quasiconformal with respect to $g$ (see \cite[Definition 4.1]{FW-Plateau-Douglas} and the explanation thereafter). It follows from \cite[Theorem 3.6]{LW20} that $u$ is a local homeomorphism. Monotonicity of $u$ implies then that $u$ is injective. Hence, $u$ is a homeomorphism as it is a continuous bijection on a compact set $M$. Using analogous statements to Theorem 2.5 and Proposition 3.5 in \cite{LW20} for the domain $(M,g)$ instead of $(\overline{D},g_{\mathrm{Eucl}})$, one can argue as in the proof of \cite[Theorem 6.1]{LW20} to obtain that $u$ is a quasisymmetric homeomorphism with respect to $g$. Note that the analogue to \cite[Theorem 2.5]{LW20} follows since $M$ admits a $(1,2)$-Poincaré inequality and is thus a Loewner space, see \cite[Theorem 9.10]{Hei01}.

It remains to show uniqueness of $(u,g)$ up to precomposition with conformal diffeomorphisms. Let $(u,g)$, $(v,h)$ be energy minimizing pairs in $\Lambda_{\text{metr}}(M,X)$. We claim that the map $\varphi:=v^{-1}\circ u\colon (M,g)\to(M,h)$ is then a conformal diffeomorphism. Indeed, for any choice of conformal charts $\psi\colon \overline{U}\to \overline{D}$ of $(M,g)$ and $\phi\colon\overline{V}\to\overline{D}$ of $(M,h)$, we can argue as in the last paragraph in the proof of \cite[Theorem 6.1]{LW20} that the transition maps $$\phi\circ v^{-1}\circ u\circ\psi^{-1}\colon \overline{D}\to\overline{D}$$ are conformal diffeomorphisms, which implies the respective property for the mapping $\varphi$. The proof of the theorem is complete.
\end{proof}

\section{Appendix}

\begin{proof}[Proof of Lemma~\ref{lem:QuotientRegularLLC}] 
Let $z\in\widehat{XY}$ and $r>0$ be arbitrary. By symmetry, we may assume $z\in X$. Observe that there exists $y\in Y$ such that $B_{\widehat{XY}}(z,r)$ is contained in $(B_{\widehat{XY}}(z,r)\cap X)\cup (B_{\widehat{XY}}(y,2r)\cap Y)$. The Ahlfors 2-regularity of $(\widehat{XY},\widehat{d})$ now follows from Lemma~\ref{lem:identity-L-Lip} and the Ahlfors 2-regularity of $X$ and $Y$.

It remains to prove that $(\widehat{XY}, \hat{d})$ is LLC. Both $X$ and $Y$ are quasiconvex (see \cite[Theorem~B.6]{Sem96}) with constants $C_X$ and $C_Y$ depending only on the LLC and Ahlfors 2-regularity constants of $X$ and $Y$, respectively. Hence, the space $(\widehat{XY}, \hat{d})$ is quasiconvex with constant $\widehat{C}:=\max\{C_X,C_Y\}$ implying that the first LLC condition holds with constant $\widehat{C}$.

Denote by $\lambda_X$ and $\lambda_Y$ the LLC-constants of $X$ and $Y$, respectively, and choose $$\hat{\lambda}\geq\max\{2,\lambda_X,\lambda_Y\}$$ such that $2\mathrm{diam}_{\hat{d}}(\widehat{XY})/\hat{\lambda}<\mathrm{diam}_{\hat{d}}(X\cap Y)$. Let $x,y\in \widehat{XY}\setminus B_{\widehat{XY}}(z,r)$. We want to prove the existence of a uniform $\lambda\geq 1$ such that $x,y$ can be joined by a continuum in $ \widehat{XY}\setminus B_{\widehat{XY}}(z,r/\lambda)$. If $x,y\in X$ or $x,y\in Y$, the statement follows from the LLC-property of $X$ or $Y$ and Lemma~\ref{lem:identity-L-Lip}. 
Consider $x\in X$, $y\in Y\setminus X$ and assume for the moment that $B:= B_{\widehat{XY}}\left(z,r/(2L\hat{\lambda}^2)\right)\subset X$. Choose any point $a\in (X\cap Y)\setminus B_{\widehat{XY}}(z,r/\hat{\lambda})$. Then there exists a continuum in $$X\setminus B_{X}\left(z,r/\hat{\lambda}^2\right)\,\subset\,\widehat{XY}\setminus B_{\widehat{XY}}\left(z,r/(L\hat{\lambda}^2)\right)\,\subset\,\widehat{XY}\setminus B$$ connecting $x$ with $a$, which can be concatenated with any continuum in $Y$ connecting $a$ with $y$ to obtain a desired path between $x$ and $y$ in $\widehat{XY}\setminus B$. If the intersection of $B_{\widehat{XY}}\left(z,r/(2L\hat{\lambda}^2)\right)$ with $Y$ is not empty, choose a point $b\in X\cap Y\cap B$ and define
$$d:=\hat{d}(b,z)<\frac{r}{2L\hat{\lambda}^2}<\frac{r}{\hat{\lambda}}.$$
It then follows from the triangle inequality that
$$d_X(b,x)\geq \hat{d}(b,x)\geq r-d \geq r-\frac{r}{\hat{\lambda}}\geq \frac{r}{\hat{\lambda}}$$
and similarly, that $d_Y(b,y)\geq r/\hat{\lambda}$. After picking a point $a\in (X\cap Y)\setminus B_{\widehat{XY}}(b, r/\hat{\lambda})$, we have the existence of continua $E\subset X\setminus B_{X}(b,r/\hat{\lambda}^2)$ connecting $x$ with $a$ respectively $F\subset Y\setminus B_{Y}(b, r/\hat{\lambda}^2)$ joining $a$ with $y$; and therefore a continuum in $$\widehat{XY}\setminus B_{\widehat{XY}}\left(b,r/(L\hat{\lambda}^2)\right)\,\subset\, \widehat{XY}\setminus B_{\widehat{XY}}\left(z, r/(L\hat{\lambda}^2)-d\right)\,\subset\, \widehat{XY}\setminus B$$ between $x$ and $y$. We thus have proven that the space $(\widehat{XY}, \hat{d})$ is LLC with constant $\lambda:=\max\{2L\hat{\lambda}^2,\widehat{C}\}.$
\end{proof}

\def\cprime{$'$} \def\cprime{$'$} \def\cprime{$'$}

\end{document}